\documentclass[12pt,oneside,reqno]{amsart}
\usepackage{mathrsfs}
\usepackage{graphics}
\usepackage{enumerate, amssymb}
\newcommand{\dif}{\mathrm{d}}

\newcommand{\be}{\begin{eqnarray}}
\newcommand{\ee}{\end{eqnarray}}
\newcommand{\ce}{\begin{eqnarray*}}
\newcommand{\de}{\end{eqnarray*}}
\newtheorem{theorem}{Theorem}[section]
\newtheorem{lemma}[theorem]{Lemma}
\newtheorem{remark}[theorem]{Remark}
\newtheorem{definition}[theorem]{Definition}
\newtheorem{proposition}[theorem]{Proposition}
\newtheorem{Example}[theorem]{Example}
\newtheorem{corollary}[theorem]{Corollary}

\def\[{{\Big[}}
\def\]{{\Big]}}
\def\<{{\langle}}
\def\>{{\rangle}}
\def\({{\Big(}}
\def\){{\Big)}}

\def\no{\nonumber}
\def\bt{\begin{theorem}}
\def\et{\end{theorem}}
\def\bl{\begin{lemma}}
\def\el{\end{lemma}}
\def\br{\begin{remark}}
\def\er{\end{remark}}
\def\bx{\begin{Example}}
\def\ex{\end{Example}}
\def\bd{\begin{definition}}
\def\ed{\end{definition}}
\def\bp{\begin{proposition}}
\def\ep{\end{proposition}}
\def\bc{\begin{corollary}}
\def\ec{\end{corollary}}

\def\cB{{\mathcal B}}
\def\cC{{\mathcal C}}
\def\cD{{\mathcal D}}
\def\cE{{\mathcal E}}
\def\cF{{\mathcal F}}
\def\cG{{\mathcal G}}

\def\cL{{\mathcal L}}
\def\cM{{\mathcal M}}
\def\cN{{\mathcal N}}

\def\cP{{\mathcal P}}

\def\mE{{\mathbb E}}

\def\mP{{\mathbb P}}
\def\mQ{{\mathbb Q}}
\def\mR{{\mathbb R}}
\def\mS{{\mathbb S}}

\def\mU{{\mathbb U}}

\def\leq{\leqslant}

\begin{document}

\allowdisplaybreaks

\title{Uniqueness for Measure-Valued Equations of Nonlinear Filtering for Stochastic Dynamical Systems with L\'evy Noises*}

\author{Huijie Qiao}

\thanks{{\it AMS Subject Classification(2010):} 60G35, 60G51, 60H10.}

\thanks{{\it Keywords:} The Zakai equation, the Kushner-Stratonovich equation, pathwise
uniqueness, uniqueness in joint law.}

\thanks{*This work was partly supported by NSF of China (No. 11001051, 11371352) and National Statistical Scientific Research Project of China (No. 2017LZ29).}

\subjclass{}

\date{}

\dedicatory{Department of Mathematics,
Southeast University\\
Nanjing, Jiangsu 211189,  China\\
hjqiaogean@seu.edu.cn}

\begin{abstract}
In the article, Zakai and Kushner-Stratonovich equations of the nonlinear filtering problem for a non-Gaussian
signal-observation system are considered. Moreover, we prove that under some general assumption, the Zakai equation has pathwise
uniqueness and uniqueness in joint law, and the Kushner-Stratonovich equation is unique in joint law.
\end{abstract}

\maketitle \rm

\section{Introduction}

Given a complete filtered probability space $(\Omega, \mathscr{F}, \{\mathscr{F}_t\}_{t\in[0,T]},\mP)$, with $T>0$ a fixed time. $X, Y$ are two processes defined on it. And as usual $X$ is difficult to observe and is called the signal process, and $Y$ is easy to observe and is called the observation process. Moreover, $Y$ contains the information about $X$. Thus, the nonlinear filtering  problem means to estimate the state of $X$ by $Y$. Precisely speaking, it is to evaluate the `filter' $\mE[F(X_t)|\mathscr{F}_t^Y]$, where $\mathscr{F}_t^Y$ is the $\sigma$-algebra generated by $\{Y_s, 0\leq s\leq t\}$ and $\mE|F(X_t)|<\infty$ for $t\in[0,T]$.

The nonlinear filtering  problem is closely related with two measure-valued equations--the Zakai and Kushner-Stratonovich
equations. Moreover, it can be said to be completely solved if solutions for the Zakai and Kushner-Stratonovich equations have uniqueness.
Therefore, in order to solve the nonlinear filtering  problem completely, studying uniqueness of the two equations is necessary.
This has been done by various authors using essentially two types of techniques.

One approach is via filtered martingale problems. In \cite{ko} Kurtz-Ocone used the technique to
prove uniqueness of solutions for the Zakai and Kushner-Stratonovich equations when $X$ is a c\`adl\`ag solution of a martingale problem and $Y$ is a continuous diffusion process. If $X$ is one jump diffusion process and $Y$ is the other jump diffusion process correlated with the Wiener process and the jump process of $X$ in one-dimensional case, Ceci-Colaneri in \cite{cckc1, cckc2} showed uniqueness of solutions for the two equations by the approach. Later the author and Duan in \cite{qd} applied the technique to uniqueness of solutions for the two equations when $X, Y$ are both multi-dimensional It\^o-L\'evy diffusion processes. The other approach is using operator techniques.
Szpirglas \cite{sz} looked like the Zakai and Kushner-Stratonovich equations as two stochastic differential equations and applied
the evolution equation of an operator to the uniqueness problem. There $X$ is a Markov process independent of the Wiener process in $Y$. In \cite{lh} Lucic-Heunis developed the technique
and studied the uniqueness problem when $X, Y$ are both continuous diffusion processes, and $X$ depends on the Wiener process in $Y$.

In the paper, we add jumps to the observation process and consider uniqueness of the two equations by means of operator techniques. In the concrete, we define weak solutions, pathwise uniqueness and uniqueness in joint law for the two measure-valued stochastic differential equations. And then, under some general assumption conditions, we prove that the Zakai equation has pathwise
uniqueness and uniqueness in joint law, and the weak solution of the Kushner-Stratonovich equation is unique in joint law. Here the signal process doesn't contain a pure jump diffusion process, since the infinitesimal generator of a pure jump diffusion process doesn't have the property that is used in the proof of a main result. Besides, we don't consider that $X$ depends on the Wiener process or the pure jump diffusion process  in $Y$ or on $Y$. And this is our future work.

It is worthwhile to mention that in \cite{qd} the author and Duan proved pathwise uniqueness for the Zakai equation. There the property is defined for strong solutions. Therefore, those assumption conditions are stronger than that here. 

This paper is arranged as follows. In Section \ref{pre}, we introduce a nonlinear filtering problem
for a It\^o-L\'evy signal-observation system and  the evolution equation for an operator.
In Section \ref{uniquezakai}, pathwise uniqueness and uniqueness in joint law of the Zakai equation
are proved. Uniqueness in joint law of the Kushner-Stratonovich equation is placed in Section \ref{uniqueks}.

\section{Preliminary}\label{pre}

In the section, we introduce some notation, terminology, concepts and known results used in the sequel.

\subsection{Notation and terminology}
Let $\cB(\mR^n)$ be the collection of all uniformly bounded Borel-measurable real-valued functions on $\mR^n$. Let $\cC_b(\mR^n)$ denote the set of all uniformly bounded continuous real-valued functions on $\mR^n$. Let $\hat{\cC}(\mR^n)$ be the collection of all members of $\cC_b(\mR^n)$ which vanish at infinity. $\cC_b^1(\mR^n)$ stands for the collection of all real-valued functions on $\mR^n$ which themself and their one-order derivatives are uniformly bounded. $\cC_c^\infty(\mR^n)$ is the collection of all real-valued functions on $\mR^n$ with continuous derivatives of all orders and compact support.

Let $\cM(\mR^n)$ be all positive bounded measures on $\mR^n$ and $\cP(\mR^n)$ be all probability measures
on $\mR^n$. For $\mu\in\cM(\mR^n)$ and a $\mathscr{B}(\mR^n)$-measurable and $\mu$-integrable function $\phi: \mR^n\rightarrow\mR$, $\mu(\phi):=\int_{\mR^n}\phi(x)\mu(\dif x)$.

\bd\label{bpclo}
Suppose that $E$ is a separable metric space. $b.p.-\lim\limits_{n\rightarrow\infty}\varphi_n=\varphi$ indicates that for $\varphi, \varphi_n\in \cB(E)$,
$\sup\limits_n\|\varphi_n\|<\infty$ (where $\|\cdot\|$ denotes sup norm in $\cB(E)$) and $\lim\limits_{n\rightarrow\infty}\varphi_n(x)=\varphi(x),\forall x\in E$.
A set $M\subset \cB(E)$ is called b.p.-closed when $b.p.-\lim\limits_{n\rightarrow\infty}\varphi_n=\varphi$
for some sequence $\{\varphi_n\}\subset M$ implies $\varphi\in M$. The b.p.-closure of $M\subset \cB(E)$
is defined to be the intersection of all b.p.-closed $M_i\subset \cB(E)$ such that $M\subset M_i$.
\ed

\subsection{A nonlinear filtering problem}\label{nonfilter}

In this subsection, we observe the nonlinear filtering problem for a non-Gaussian signal-observation system,  and
state the Zakai and Kushner-Stratonovich equations. (c.f.\cite{qd})

Fix $T>0$, and consider the following signal-observation $(X_t, Y_t)$  system   on $\mR^n\times\mR^m$:
\be\left\{\begin{array}{l}
\dif X_t=b_1(t,X_t)\dif t+\sigma_1(t,X_t)\dif B_t, \\
\dif Y_t=b_2(t,X_t)\dif t+\sigma_2(t)\dif W_t+\int_{\mU_0}f_2(t,u)\tilde{N}_{\lambda}(\dif t, \dif u)
 +\int_{\mU\setminus\mU_0}g_2(t,u)N_{\lambda}(\dif t, \dif u),
\end{array}
\right. \no\\
0\leq t\leq T,\qquad \label{Eq1}
\ee
where $B, W$ are $d$-dimensional and $m$-dimensional Brownian motions, respectively, and $N_{\lambda}(\dif t, \dif u)$ is an
integer-valued random measure with a predictable compensator $\lambda(t,X_{t-},u)\dif t\nu(\dif u)$. Here the
function $\lambda(t,x,u)\in(0,1)$, and $\nu$ is a $\sigma$-finite measure defined on a measurable
space ($\mU, \mathscr{U}$) with $\nu(\mU\setminus\mU_0)<\infty$ for
$\mU_0\in\mathscr{U}$. That is, $\tilde{N}_\lambda(\dif t,\dif u):=N_\lambda(\dif t,\dif u)-\lambda(t,X_{t-},u)\dif t\nu(\dif u)$ is its
compensated martingale measure. Moreover, $B_t, W_t, N_{\lambda}$ are mutually independent. The initial
value $X_0$ is assumed to be a random variable independent of $Y_0, B_t, W_t, N_{\lambda}$.

The mappings $b_1: [0,T]\times\mR^n\mapsto\mR^n$, $b_2: [0,T]\times\mR^n\mapsto\mR^m$, $\sigma_1: [0,T]\times\mR^n\mapsto\mR^{n\times d}$,
$\sigma_2: [0,T]\mapsto\mR^{m\times m}$, $f_2: [0,T]\times\mU_0\mapsto\mR^m$ and $g_2: [0,T]\times(\mU\setminus\mU_0)\mapsto\mR^m$
are all Borel measurable. We make the following  assumptions.

{\bf Assumption 1.}

\begin{enumerate}[($\mathbf{H}^1_{b_1, \sigma_1}$)]
\item  For $t\in[0,T]$, $b_1(t,x)$ is continuous in $x$, and $\sigma_1(t,x)$ is bi-continuous in $(t, x)$.
\end{enumerate}

\begin{enumerate}[($\mathbf{H}^2_{b_1, \sigma_1}$)]
\item  There exists a non-negative constant $L_1$ such that
\ce
|b_1(t, x)|^2+\|\sigma_1(t, x)\|^2\leq L_1(1+|x|)^2, \quad t\in[0,T], x\in\mR^n.
\de
\end{enumerate}

\begin{enumerate}[($\mathbf{H}^1_{b_2, \sigma_2, f_2}$)]
\item  $\sigma_2(t)$ is invertible for $t\in[0,T]$, $b_2, \sigma_2, \sigma^{-1}_2$ are bounded by a positive constant $L_2$, and
\ce
\int_0^T\int_{\mU_0}|f_2(s,u)|^2\nu(\dif u)\dif s<\infty.
\de
\end{enumerate}

\medskip

By \cite[Chapter 5, Theorem 175]{si}, the system (\ref{Eq1}) has a weak solution denoted
by $(X_t,Y_t)$. Set
\ce
\Lambda^{-1}_t:&=&\exp\bigg\{-\int_0^t\left(\sigma_2^{-1}(s)b_2(s,X_s)\right)^i\dif W^i_s-\frac{1}{2}\int_0^t
\left|\sigma_2^{-1}(s)b_2(s,X_s)\right|^2\dif s\\
&&\quad\qquad -\int_0^t\int_{\mU_0}\log\lambda(s,X_{s-},u)\tilde{N}_{\lambda}(\dif s, \dif u)\\
&&\quad\qquad-\int_0^t\int_{\mU_0}\left[\log\lambda(s,X_{s-},u)+\frac{(1-\lambda(s,X_{s-},u))}{\lambda(s,X_{s-},u)}\right]\lambda(s,X_{s-},u)\nu(\dif u)\dif s\bigg\}.
\de
Here and hereafter, we use the convention that repeated indices imply summation.

\medskip

{\bf Assumption 2.}
There exists a positive function $L(u)$ satisfying
\ce
\int_{\mU_0}\frac{\left(1-L(u)\right)^2}{L(u)}\nu(\dif u)<\infty
\de
such that $0<l\leq L(u)<\lambda(t,x,u)<1$ for $u\in\mU_0$, where $l$ is a constant.

\medskip

Under {\bf Assumption 2.}, it holds that
\ce
&&\mE\left[\exp\left\{\int_0^T\int_{\mU_0}\frac{\left(1-\lambda(s,X_s,u)\right)^2}{\lambda(s,X_s,u)}\nu(\dif u)\dif s\right\}\right]\\
&<&\exp\left\{\int_0^T\int_{\mU_0}\frac{\left(1-L(u)\right)^2}{L(u)}\nu(\dif u)\dif s\right\}\\
&<&\infty.
\de
Thus, by the similar deduction to that in \cite{qd}, we know that $\Lambda^{-1}_t$ is an exponential martingale. Define a probability measure $\tilde{\mP}$ via
$$
\frac{\dif \tilde{\mP}}{\dif \mP}=\Lambda^{-1}_T.
$$
By the Girsanov theorem for Brownian motions and random measures, under the measure
$\tilde{\mP}$ the system (\ref{Eq1}) is transformed as
\be\left\{\begin{array}{l}
\dif X_t=b_1(t, X_t)\dif t+\sigma_1(t, X_t)\dif B_t,\\
\dif Y_t=\sigma_2(t)\dif \tilde{W}_t+\int_{\mU_0}f_2(t,u)\tilde{N}(\dif t, \dif u)+\int_{\mU\setminus\mU_0}g_2(t,u)N_{\lambda}(\dif t, \dif u),
\end{array}
\right. \label{Eq2}
\ee
where
\ce
\tilde{W}_t:=W_t+\int_0^t\sigma_2^{-1}(s)b_2(s,X_s)\dif s, \quad
\tilde{N}(\dif t, \dif u):=N_\lambda(\dif t, \dif u)-\dif t\nu(\dif u).
\de
Moreover, under the measure $\tilde{\mP}$, $\tilde{W}$ is a Brownian motion and $\tilde{N}$ is a Poisson compensated martingale measure.

Set
\ce
\tilde{\mP}_t(F):=\tilde{\mE}[F(X_t)\Lambda_t|\mathscr{F}_t^Y], \quad F\in\cB(\mR^n), t\in[0,T],
\de
where $\tilde{\mE}$ denotes expectation under the measure $\tilde{\mP}$ and $\mathscr{F}_t^Y$ is the $\sigma$-algebra
generated by $\{Y_s, 0\leq s\leq t\}$. The equation satisfied by $\tilde{\mP}_t(F)$
is called the Zakai equation. Based on Theorem 3.1 in \cite{qd}, we have the following result.

\bt (Zakai equation) \label{zakait}
For $F\in\cC_c^\infty(\mR^n)\cup\{1\}$, the Zakai equation of the system (\ref{Eq1}) is given by
\be
\tilde{\mP}_t(F)&=&\tilde{\mP}_0(F)+\int_0^t\tilde{\mP}_s(\cL_s F)\dif s+\int_0^t\tilde{\mP}_s\left(F\left(\sigma_2^{-1}(s)b_2(s,\cdot)\right)^i\right)\dif\tilde{W}^i_s\no\\
&&+\int_0^t\int_{\mU_0}\tilde{\mP}_{s-}\left(F(\lambda(s,\cdot,u)-1)\right)\tilde{N}(\dif s, \dif u),
\label{zakaieq1}
\ee
where $\cL_t$ is the infinitesimal generator of $X_t$ and is given by
\ce
(\cL_t F)(x)&:=&\frac{\partial F(x)}{\partial x_i}b^i_1(t,x)+\frac{1}{2}\frac{\partial^2F(x)}{\partial x_i\partial x_j}
\sigma_1^{ik}(t, x) \sigma_1^{jk}(t, x).
\de
\et

Besides, set
\ce
\mP_t(F):=\mE[F(X_t)|\mathscr{F}_t^Y], \quad F\in\cB(\mR^n), t\in[0,T],
\de
and then it follows from the Kallianpur-Striebel formula that
\be
\mP_t(F)=\mE[F(X_t)|\mathscr{F}_t^Y]=\frac{\tilde{\mE}[F(X_t)\Lambda_t|\mathscr{F}_t^Y]}
{\tilde{\mE}[\Lambda_t|\mathscr{F}_t^Y]}=\frac{\tilde{\mP}_t(F)}{\tilde{\mP}_t(1)}.
\label{bayf}
\ee
By \cite[Theorem 3.2]{qd}, we obtain the following Kushner-Stratonovich equation satisfied by $\mP_t(F)$.

\bt  (Kushner-Stratonovich equation)  \label{ks}
For $F\in\cC_c^\infty(\mR^n)$, $\mP_t(F)$ solves the following equation
\be
\mP_t(F)&=&\mP_0(F)+\int_0^t\mP_s(\cL_s F)\dif s\no\\
&&+\int_0^t\bigg(\mP_s\left(F\left(\sigma_2^{-1}(s)b_2(s,\cdot)\right)^i\right)
-\mP_s\left(F\right)\mP_s\left(\left(\sigma_2^{-1}(s)b_2(s,\cdot)\right)^i\right)\bigg)\dif \bar{W}^i_s\no\\
&&+\int_0^t\int_{\mU_0}\frac{\mP_{s-}\left(F\lambda(s,\cdot,u)\right)-\mP_{s-}\left(F\right)\mP_{s-}
\left(\lambda(s,\cdot,u)\right)}{\mP_{s-}\left(\lambda(s,\cdot,u)\right)}\tilde{\bar{N}}(\dif s, \dif u),
\label{kseq}
\ee
where $\bar{W}_t:=\tilde{W}_t-\int_0^t\mP_s\left(\sigma_2^{-1}(s)b_2(s,\cdot)\right)\dif s$ and
$\tilde{\bar{N}}(\dif t, \dif u)=N_\lambda(\dif t, \dif u)-\mP_{t}\left(\lambda(t,\cdot,u)\right)\nu(\dif u)\dif t$.
\et

By Theorem VI.8.4 in \cite{rw}, $\bar{W}_t$ is a $\{\mathscr{F}_t^Y\}$-Brownian motion under $\mP$. Based on
the tower property of conditional expectation, we know that $\tilde{\bar{N}}(\dif t, \dif u)$ is the compensated martingale measure
for the random measure $N_\lambda(\dif t, \dif u)$ with the predictable compensator $\mP_{t}\left(\lambda(t,\cdot,u)\right)\nu(\dif u)\dif t$
under $\mP$.

\subsection{An evolution equation for an operator}\label{evoequ}

In the subsection, we introduce an evolution equation for an operator and prove a related result used in the following section.

Suppose that $E$ is a complete separable metric space, and an operator $\cL$ is defined on $\cB(E)$ with domain $\cD(\cL)$. If there exists a family $\{\mu_t, t\in[0,T]\}$ such that (i) $\mu_t\in\cM(E)$ for any $t\in[0,T]$ and
$\mu_0\in\cP(E)$; (ii) for $B\in\mathscr{B}(E)$, $\mu_t(B)$ is Borel measurable in $t$; (iii) for any
$\varphi\in\cD(\cL)$, it holds that $\int_0^t|\mu_s(\cL\varphi)|\dif s<\infty$ for any $t\in[0,T]$ and
$$
\mu_t(\varphi)=\mu_0(\varphi)+\int_0^t\mu_s(\cL\varphi)\dif s, \quad t\in[0,T];
$$
we call $\{\mu_t, t\in[0,T]\}$ a $\cM(E)$-valued solution of the evolution equation for $(\cL, \cD(\cL))$. And
uniqueness in the class of $\cM(E)$-valued solutions over the interval $[0,T]$, of the evolution equation for $(\cL, \cD(\cL))$,
means that if there exist two such solutions $\{\mu^1_t, t\in[0,T]\}$ and $\{\mu^2_t, t\in[0,T]\}$ with $\mu^1_0=\mu^2_0$,
then $\mu^1_t=\mu^2_t$ for any $t\in[0,T]$. In the following, to a type of special operators, we give some conditions to justify
uniqueness in the class of $\cM(E)$-valued solutions over the interval $[0,T]$ of the evolution equations for them.

Define an operator $\mathscr{L}_t$ by
\ce
\cD(\mathscr{L}_t):=span\{1,\cC_c^\infty(\mR^{2n})\},
\de
\ce
\mathscr{L}_t\phi(x):=\frac{\partial\phi(x)}{\partial x_i}\beta^i(t,x)+\frac{1}{2}\frac{\partial^2\phi(x)}{\partial x_i\partial x_j}\alpha^{ij}(t,x),
\quad t\in[0,T], x\in\mR^{2n}, \phi\in\cD(\mathscr{L}_t),
\de
where $\beta: [0,T]\times\mR^{2n}\mapsto\mR^{2n}$ is Borel measurable and $\alpha: [0,T]\times\mR^{2n}\mapsto\mS^{2n}_{+}$ is continuous.
Here $\mS^{2n}_{+}$ denotes the class of all nonnegative definite $2n\times 2n$ matrices.

\bt\label{unmeev}
Suppose that there exists a positive constant $L_3$ such that for $(t,x)\in[0,T]\times\mR^{2n}$,
\ce
|\beta(t,x)|\leq L_3(1+|x|), \quad \|\alpha(t,x)\|\leq L_3(1+|x|^2).
\de
If $\gamma: [0,T]\times\mR^{2n}\mapsto\mR$ is a uniformly bounded Borel measurable function, the evolution equation for $(\mathscr{L}_t-\gamma(t,\cdot), \cD(\mathscr{L}_t))$ is unique in the class of $\cM(\mR^{2n})$-valued solutions over the interval $[0,T]$.
\et
\begin{proof}
Set $E:=[0,T]\times\mR^{2n}$ and $\cD':=span\{h\phi; h\in \cC^1_b([0,T]), \phi\in\cD(\mathscr{L}_t)\}$,
and then $\cD'\subset \cC_b(E)$. Define an operator $\mathscr{L}$ by
\ce
\cD(\mathscr{L}):=\cD',
\de
\ce
\mathscr{L}\varphi(t,x):=h'(t)\phi(x)+h(t)\mathscr{L}_t\phi(x), \quad \varphi(t,x)=h(t)\phi(x)\in\cD(\mathscr{L}).
\de
Thus, $\{\rho_t, t\in[0,T]\}$ is a $\cM(\mR^{2n})$-valued solution of the evolution equation for $(\mathscr{L}_t-\gamma(t,\cdot), \cD(\mathscr{L}_t))$
if and only if $\{\delta_t\times\rho_t, t\in[0,T]\}$ is a $\cM(E)$-valued solution of the evolution equation for $(\mathscr{L}-\gamma(\cdot,\cdot), \cD(\mathscr{L}))$, where $\delta_t$ is the Dirac measure at $t$. In fact, if $\{\rho_t, t\in[0,T]\}$ is a $\cM(\mR^{2n})$-valued solution of the evolution equation for $(\mathscr{L}_t-\gamma(t,\cdot), \cD(\mathscr{L}_t))$, by the above definition, it holds that for any $\phi\in\cD(\mathscr{L}_t)$,
$$
\rho_t(\phi)=\rho_0(\phi)+\int_0^t\rho_s((\mathscr{L}_s-\gamma(s,\cdot))\phi)\dif s, \quad t\in[0,T].
$$
Note that for any $h\in \cC^1_b([0,T])$
$$
\delta_t(h)=\delta_0(h)+\int_0^t\delta_s(h')\dif s, \quad t\in[0,T].
$$
Combining the two equalities with integration by parts, we obtain that
$$
\delta_t(h)\rho_t(\phi)=\delta_0(h)\rho_0(\phi)+\int_0^t\delta_s(h')\rho_s(\phi)\dif s+\int_0^t\delta_s(h)\rho_s((\mathscr{L}_s
-\gamma(s,\cdot))\phi)\dif s,
$$
and
$$
\delta_t\rho_t(\varphi)=\delta_0\rho_0(\varphi)+\int_0^t\delta_s\rho_s((\mathscr{L}-\gamma(\cdot,\cdot))\varphi)\dif s, \quad \varphi(t,x)=h(t)\phi(x)\in\cD(\mathscr{L}).
$$
Thus, $\{\delta_t\times\rho_t, t\in[0,T]\}$ is a $\cM(E)$-valued solution of the evolution equation for $(\mathscr{L}-\gamma(\cdot,\cdot), \cD(\mathscr{L}))$. Conversely, it is simple by taking $h(t)=1$ for $t\in[0,T]$.

Since $\alpha, \beta, \gamma$ satisfy these assumptions, by Theorem 3.6 in \cite{lh}, the evolution equation for $(\mathscr{L}-\gamma(\cdot,\cdot), \cD(\mathscr{L}))$ has uniqueness in the class of $\cM(E)$-valued solutions over the
interval $[0,T]$. From this, we know that the evolution equation for $(\mathscr{L}_t-\gamma(t,\cdot), \cD(\mathscr{L}_t))$ is unique in the class of $\cM(\mR^{2n})$-valued solutions over the interval $[0,T]$. The proof is completed.
\end{proof}

\br
In the proof of the above theorem, Theorem 3.6 in \cite{lh} is used. Certainly speaking, we apply its modified version to $\alpha, \beta, \gamma$,
since $\alpha$ is required to be strictly positive definite there. After checking the proof of Theorem 3.6 carefully, we find that nonnegative definite
property of $\alpha$ is only used. That is, when $\alpha$ is nonnegative definite, the result in Theorem 3.6 is still right.
\er

\section{Pathwise uniqueness and uniqueness in joint law of the Zakai equation}\label{uniquezakai}

In the section, we define weak solutions, pathwise uniqueness and uniqueness in joint law of the Zakai equation in Subsection \ref{nonfilter}. And then, we prove that the Zakai equation has pathwise
uniqueness and uniqueness in joint law. Let us start with some notations.

\bd\label{soluzakai}
If there exists the pair $\{(\hat{\Omega}, \hat{\mathscr{F}}, \{\hat{\mathscr{F}}_t\}_{t\in[0,T]}, \hat{\mP}), (\hat{\mu}_t,
\hat{W}_t, \hat{N}(\dif t, \dif u))\}$ such that the following hold:

(i) $(\hat{\Omega}, \hat{\mathscr{F}}, \{\hat{\mathscr{F}}_t\}_{t\in[0,T]},\hat{\mP})$ is a complete filtered
probability space;

(ii) $\hat{\mu}_t$ is a $\cM(\mR^n)$-valued $\hat{\mathscr{F}}_t$-adapted c\`adl\`ag process and $\hat{\mu}_0\in\cP(\mR^n)$;

(iii) $\hat{W}_t$ is a $m$-dimensional $\hat{\mathscr{F}}_t$-adapted Brownian motion;

(iv) $\hat{N}(\dif t, \dif u)$ is a Poisson random measure with a predictable compensator $\dif t\nu(\dif u)$;

(v) $(\hat{\mu}_t, \hat{W}_t, \hat{N}(\dif t, \dif u))$ satisfies the following equation
\be
\hat{\mu}_t(F)&=&\hat{\mu}_0(F)+\int_0^t\hat{\mu}_s(\cL_s F)\dif s+\int_0^t\hat{\mu}_s\left(F\left(\sigma_2^{-1}(s)b_2(s,\cdot)\right)^i\right)\dif\hat{W}^i_s\no\\
&&+\int_0^t\int_{\mU_0}\hat{\mu}_{s-}\(F(\lambda(s,\cdot,u)-1)\)\tilde{\hat{N}}(\dif s, \dif u), \quad F\in \cC_c^\infty(\mR^n)\cup\{1\},
\label{zakaieq2}
\ee
where $\tilde{\hat{N}}(\dif t, \dif u):=\hat{N}(\dif t, \dif u)-\dif t\nu(\dif u)$,
then $\{(\hat{\Omega}, \hat{\mathscr{F}}, \{\hat{\mathscr{F}}_t\}_{t\in[0,T]},\hat{\mP}), (\hat{\mu}_t,
\hat{W}_t, \hat{N}(\dif t, \dif u))\}$ is called a weak solution of the Zakai equation.
\ed

By the deduction in Subsection \ref{nonfilter}, it is obvious that $\{(\Omega, \mathscr{F}, \{\mathscr{F}_t\}_{t\in[0,T]}, \tilde{\mP}),
(\tilde{\mP}_t, \tilde{W}_t, \\ N_\lambda(\dif t, \dif u))\}$ is a weak solution of the Zakai equation.

\bd\label{paunzakai}
Pathwise uniqueness of the Zakai equation means that if there exist two weak solutions $\{(\hat{\Omega}, \hat{\mathscr{F}}, \{\hat{\mathscr{F}}_t\}_{t\in[0,T]}, \hat{\mP}), (\hat{\mu}^1_t,\hat{W}_t, \hat{N}(\dif t, \dif u))\}$ and $\{(\hat{\Omega}, \hat{\mathscr{F}}, \{\hat{\mathscr{F}}_t\}_{t\in[0,T]}, \hat{\mP}), \\(\hat{\mu}^2_t,
\hat{W}_t, \hat{N}(\dif t, \dif u))\}$ with $\hat{\mP}\{\hat{\mu}^1_0=\hat{\mu}^2_0\}=1$, then
$$
\hat{\mu}^1_t=\hat{\mu}^2_t, \quad t\in[0,T], ~ a.s.\hat{\mP}.
$$
\ed

\bd\label{launzakai}
Uniqueness in joint law of the Zakai equation means that if there exist two weak solutions $\{(\hat{\Omega}^1, \hat{\mathscr{F}}^1, \{\hat{\mathscr{F}}^1_t\}_{t\in[0,T]}, \hat{\mP}^1), (\hat{\mu}^1_t,\hat{W}^1_t, \hat{N}^1(\dif t, \dif u))\}$ and $\{(\hat{\Omega}^2, \hat{\mathscr{F}}^2, \{\hat{\mathscr{F}}^2_t\}_{t\in[0,T]}, \hat{\mP}^2),\\ (\hat{\mu}^2_t,
\hat{W}^2_t, \hat{N}^2(\dif t, \dif u))\}$ with $\hat{\mP}^1\circ(\hat{\mu}^1_0)^{-1}=\hat{\mP}^2\circ(\hat{\mu}^2_0)^{-1}$, then
$\{(\hat{\mu}^1_t,\hat{W}^1_t, \hat{N}^1(\dif t, \dif u)), t\in[0,T]\}$ and $\{(\hat{\mu}^2_t,\hat{W}^2_t, \hat{N}^2(\dif t, \dif u)), t\in[0,T]\}$
have the same finite-dimensional distributions.
\ed

\medskip

Now, it is the position to state and prove the first main result in the section.

\bt\label{unipalazakai}
Suppose that {\bf Assumption 1.-2.} are satisfied. Then

(i) the Zakai equation has pathwise uniqueness,

(ii) the Zakai equation has uniqueness in joint law.
\et
{\bf Proof of Theorem \ref{unipalazakai} (i)}: Taking two weak solutions $\{(\hat{\Omega}, \hat{\mathscr{F}}, \{\hat{\mathscr{F}}_t\}_{t\in[0,T]}, \hat{\mP}), (\hat{\mu}^1_t,\hat{W}_t,\\ \hat{N}(\dif t, \dif u))\}$ and $\{(\hat{\Omega}, \hat{\mathscr{F}}, \{\hat{\mathscr{F}}_t\}_{t\in[0,T]}, \hat{\mP}), (\hat{\mu}^2_t, \hat{W}_t, \hat{N}(\dif t, \dif u))\}$ with $\hat{\mP}\{\hat{\mu}^1_0=\hat{\mu}^2_0\}=1$, we compute
$\hat{\mE}|\hat{\mu}^1_t(F)-\hat{\mu}^2_t(F)|^2$ for any $t\in[0,T]$ and $F\in \cC_c^\infty(\mR^n)$, where
$\hat{\mE}$ is the expectation under $\hat{\mP}$. By elementary calculation, it holds that
\be
\hat{\mE}|\hat{\mu}^1_t(F)-\hat{\mu}^2_t(F)|^2=\hat{\mE}\hat{\mu}^1_t(F)\hat{\mu}^1_t(F)
-2\hat{\mE}\hat{\mu}^1_t(F)\hat{\mu}^2_t(F)+\hat{\mE}\hat{\mu}^2_t(F)\hat{\mu}^2_t(F).
\label{ccs}
\ee

Next, we compute $\hat{\mE}\hat{\mu}^1_t(F)\hat{\mu}^1_t(F)$. Set for any $t\in[0,T]$
\be
&&\hat{\mu}^{11}_t(A_1\times A_2):=\hat{\mu}^1_t(A_1)\hat{\mu}^1_t(A_2), \quad A_1, A_2\in\mathscr{B}(\mR^n),\label{me1}\\
&&\bar{\hat{\mu}}^{11}_t(A_1\times A_2):=\hat{\mE}\hat{\mu}^{11}_t(A_1\times A_2)=\hat{\mE}\hat{\mu}^1_t(A_1)\hat{\mu}^1_t(A_2),\label{me2}
\ee
and then $\{\bar{\hat{\mu}}^{11}_t, t\in[0,T]\}$ is a measure family on $(\mR^{2n}, \mathscr{B}(\mR^{2n}))$. Furthermore, it
has some good properties. Firstly, note that by Lemma \ref{gj} for any $t\in[0,T]$
\ce
&&\bar{\hat{\mu}}^{11}_t(\mR^{2n})=\hat{\mE}\hat{\mu}^1_t(\mR^n)\hat{\mu}^1_t(\mR^n)=\hat{\mE}\hat{\mu}^1_t(1)\hat{\mu}^1_t(1)<\infty, \\
&&\bar{\hat{\mu}}^{11}_0(\mR^{2n})=\hat{\mE}\hat{\mu}^1_0(\mR^n)\hat{\mu}^1_0(\mR^n)=1.
\de
Thus $\bar{\hat{\mu}}^{11}_t\in\cM(\mR^{2n})$ for any $t\in[0,T]$ and $\bar{\hat{\mu}}^{11}_0\in\cP(\mR^{2n})$. Secondly, by
the Dynkin class theorem and the Fubini theorem, we know that for $\Gamma\in\mathscr{B}(\mR^{2n})$, $\bar{\hat{\mu}}^{11}_t(\Gamma)$
is Borel measurable in $t$. Thirdly, it follows from (\ref{me1}) (\ref{me2}) that
\ce
&&\hat{\mu}^{11}_t(F_1\otimes F_2)=\hat{\mu}^1_t(F_1)\hat{\mu}^1_t(F_2), \quad F_1, F_2\in\cB(\mR^n),\\
&&\bar{\hat{\mu}}^{11}_t(F_1\otimes F_2)=\hat{\mE}\hat{\mu}^{11}_t(F_1\otimes F_2)=\hat{\mE}\hat{\mu}^1_t(F_1)\hat{\mu}^1_t(F_2),
\de
where $F_1\otimes F_2$ is the tensor product of $F_1$ and $F_2$, i.e. $F_1\otimes F_2(x_1, x_2)=F_1(x_1)F_2(x_2)$ for $x_1, x_2\in\mR^n$.
So, $\hat{\mE}\hat{\mu}^1_t(F)\hat{\mu}^1_t(F)=\bar{\hat{\mu}}^{11}_t(F\otimes F)$.

Let us deal with $\hat{\mE}\hat{\mu}^1_t(F)\hat{\mu}^2_t(F), \hat{\mE}\hat{\mu}^2_t(F)\hat{\mu}^2_t(F)$. By the same way as $\bar{\hat{\mu}}^{11}_t$, one could define $\bar{\hat{\mu}}^{12}_t, \bar{\hat{\mu}}^{22}_t$. Thus, $\hat{\mE}\hat{\mu}^1_t(F)\hat{\mu}^2_t(F)=\bar{\hat{\mu}}^{12}_t(F\otimes F), \hat{\mE}\hat{\mu}^2_t(F)\hat{\mu}^2_t(F)=\bar{\hat{\mu}}^{22}_t(F\otimes F)$.
And then (\ref{ccs}) is written as
\be
\hat{\mE}|\hat{\mu}^1_t(F)-\hat{\mu}^2_t(F)|^2=\bar{\hat{\mu}}^{11}_t(F\otimes F)-2\bar{\hat{\mu}}^{12}_t(F\otimes F)+\bar{\hat{\mu}}^{22}_t(F\otimes F).
\label{ccs2}
\ee

In the following, we compare $\bar{\hat{\mu}}^{11}_t, \bar{\hat{\mu}}^{12}_t, \bar{\hat{\mu}}^{22}_t$. First of all, consider $\bar{\hat{\mu}}^{11}_t$.
Note that for $F_1, F_2\in \cC_c^\infty(\mR^n)$, $\hat{\mu}^1_t(F_i), i=1, 2,$ satisfy Eq.(\ref{zakaieq2}). Applying the It\^o formula to $\hat{\mu}^1_t(F_1)\hat{\mu}^1_t(F_2)$ and taking the expectation on two sides, we obtain that
\ce
\hat{\mE}\hat{\mu}^1_t(F_1)\hat{\mu}^1_t(F_2)&=&\hat{\mE}\hat{\mu}^1_0(F_1)\hat{\mu}^1_0(F_2)+\int_0^t\hat{\mE}\hat{\mu}^1_s(\cL_s F_1)\hat{\mu}^1_s(F_2)\dif s
+\int_0^t\hat{\mE}\hat{\mu}^1_s(F_1)\hat{\mu}^1_s(\cL_s F_2)\dif s\\
&&+\int_0^t\hat{\mE}\hat{\mu}^1_s\left(F_1\left(\sigma_2^{-1}(s)b_2(s,\cdot)\right)^i\right)
\hat{\mu}^1_s\left(F_2\left(\sigma_2^{-1}(s)b_2(s,\cdot)\right)^i\right)\dif s\\
&&+\int_0^t\int_{\mU_0}\hat{\mE}\hat{\mu}^1_{s}\(F_1(\lambda(s,\cdot,u)-1)\)\hat{\mu}^1_{s}\(F_2(\lambda(s,\cdot,u)-1)\)
\nu(\dif u)\dif s,
\de
and
\ce
\bar{\hat{\mu}}^{11}_t(F_1\otimes F_2)&=&\bar{\hat{\mu}}^{11}_0(F_1\otimes F_2)+\int_0^t\bar{\hat{\mu}}^{11}_s(\cL_s F_1\otimes F_2)\dif s+\int_0^t\bar{\hat{\mu}}^{11}_s(F_1\otimes \cL_s F_2)\dif s\\
&&+\int_0^t\bar{\hat{\mu}}^{11}_s\left(F_1\left(\sigma_2^{-1}(s)b_2(s,\cdot)\right)^i\otimes F_2\left(\sigma_2^{-1}(s)b_2(s,\cdot)\right)^i\right)\dif s\\
&&+\int_0^t\int_{\mU_0}\bar{\hat{\mu}}^{11}_s\(F_1(\lambda(s,\cdot,u)-1)\otimes F_2(\lambda(s,\cdot,u)-1)\)\nu(\dif u)\dif s\\
&=&\bar{\hat{\mu}}^{11}_0(F_1\otimes F_2)+\int_0^t\bar{\hat{\mu}}^{11}_s\left(\cL_s F_1\otimes F_2+F_1\otimes \cL_s F_2\right)\dif s\\
&&+\int_0^t\bar{\hat{\mu}}^{11}_s\left(F_1\left(\sigma_2^{-1}(s)b_2(s,\cdot)\right)^i\otimes F_2\left(\sigma_2^{-1}(s)b_2(s,\cdot)\right)^i\right)\dif s\\
&&+\int_0^t\bar{\hat{\mu}}^{11}_s\left(\int_{\mU_0}F_1(\lambda(s,\cdot,u)-1)\otimes F_2(\lambda(s,\cdot,u)-1)\nu(\dif u)\right)\dif s,
\de
where the expression of the last term in the above equality is based on {\bf Assumption 2.}. Set
\ce
\bar{\mathscr{L}}_t(F_1\otimes F_2)&:=&\cL_t F_1\otimes F_2+F_1\otimes \cL_t F_2,\\
-\bar{\gamma}(t, \cdot)(F_1\otimes F_2)&:=&F_1\left(\sigma_2^{-1}(t)b_2(t,\cdot)\right)^i\otimes F_2\left(\sigma_2^{-1}(t)b_2(t,\cdot)\right)^i\\
&&+\int_{\mU_0}F_1(\lambda(t,\cdot,u)-1)\otimes F_2(\lambda(t,\cdot,u)-1)\nu(\dif u),
\de
and then the above equality could be written as
\be
\bar{\hat{\mu}}^{11}_t(F_1\otimes F_2)=\bar{\hat{\mu}}^{11}_0(F_1\otimes F_2)+\int_0^t\bar{\hat{\mu}}^{11}_s\left((\bar{\mathscr{L}}_s-\bar{\gamma}(s, \cdot))(F_1\otimes F_2)\right)\dif s.
\label{zhji}
\ee

Let us observe Eq.(\ref{zhji}). For $t\in[0,T], x=(x_1, x_2), x_1, x_2\in\mR^n$, set
\ce
&&\bar{a}(t, x):=\left(\begin{array}{c}
\sigma_1(t,x_1)\sigma_1^T(t,x_1) \quad\qquad\qquad\qquad 0\\
0 \qquad\qquad\qquad\quad \sigma_1(t,x_2)\sigma_1^T(t,x_2)
\end{array}
\right), \quad
\bar{b}(t,x):=\left(\begin{array}{c}
b_1(t,x_1)\\
b_1(t,x_2)
\end{array}
\right),\\
&&\bar{h}(t, x):=\left(\sigma_2^{-1}(t)b_2(t,x_1)\right)^i\left(\sigma_2^{-1}(t)b_2(t,x_2)\right)^i,\\
&&\bar{\lambda}(t, x):=\int_{\mU_0}(\lambda(t,x_1,u)-1)(\lambda(t,x_2,u)-1)\nu(\dif u),
\de
and then by some calculations $\bar{\mathscr{L}}_t, \bar{\gamma}(t, \cdot)$ could be expressed as
\ce
&&\bar{\mathscr{L}}_t(F_1\otimes F_2)(x)=\frac{\partial(F_1\otimes F_2)(x)}{\partial x^i}\bar{b}^i(t,x)+\frac{1}{2}\frac{\partial^2(F_1\otimes F_2)(x)}{\partial x^i\partial x^j}\bar{a}^{ij}(t,x),\\
&&\bar{\gamma}(t, x)(F_1\otimes F_2)(x)=\bar{h}(t, x)(F_1\otimes F_2)(x)+\bar{\lambda}(t, x)(F_1\otimes F_2)(x).
\de
And put
$$
\cG:=span\{F_1\otimes F_2; F_1, F_2\in \cC_c^\infty(\mR^n)\},
$$
and then Eq.(\ref{zhji}) is simply written as
\be
\bar{\hat{\mu}}^{11}_t(\phi)=\bar{\hat{\mu}}^{11}_0(\phi)+\int_0^t\bar{\hat{\mu}}^{11}_s((\bar{\mathscr{L}}_s-\bar{\gamma}(s,\cdot))\phi)\dif s, \quad \phi\in\cG.
\label{tenpro}
\ee
Moreover, we can expand the operator $\bar{\mathscr{L}}_t$ still denoted by $\bar{\mathscr{L}}_t$ such that
$\cD(\bar{\mathscr{L}}_t)=span\{1, \cC_c^\infty(\mR^{2n})\}$. {\bf Claim:} $\bar{\hat{\mu}}^{11}$
solves Eq.(\ref{tenpro}) for $\phi\in\cD(\bar{\mathscr{L}}_t)$.

On one side, by some similar proofs to $\bar{\hat{\mu}}^{11}$, we obtain that $\bar{\hat{\mu}}^{12}, \bar{\hat{\mu}}^{22}$ also
solve Eq.(\ref{tenpro}) for $\phi\in\cD(\bar{\mathscr{L}}_t)$. On the other side,  note that Eq.(\ref{tenpro}) for $\phi\in\cD(\bar{\mathscr{L}}_t)$
is exactly the evolution equation for $(\bar{\mathscr{L}}_t-\bar{\gamma}(t,\cdot), \cD(\bar{\mathscr{L}}_t))$. And under {\bf Assumption 1.-2.}, Theorem \ref{unmeev} admits us to get that Eq.(\ref{tenpro}) is unique in the class of $\cM(\mR^{2n})$-valued solutions over the interval $[0,T]$. So,
\be
\bar{\hat{\mu}}^{11}_t=\bar{\hat{\mu}}^{12}_t=\bar{\hat{\mu}}^{22}_t, \quad t\in[0,T].
\label{ccs3}
\ee
Combining (\ref{ccs3}) with (\ref{ccs2}), we have that $\hat{\mE}|\hat{\mu}^1_t(F)-\hat{\mu}^2_t(F)|^2=0$ and $\hat{\mu}^1_t(F)=\hat{\mu}^2_t(F)$ a.s.$\hat{\mP}$. Based on density of $\cC_c^\infty(\mR^n)$ in $\hat{\cC}(\mR^n)$ and separablility of $\hat{\cC}(\mR^n)$, we furthermore have $\hat{\mu}^1_t(F)=\hat{\mu}^2_t(F), F\in\hat{\cC}(\mR^n) $ a.s.$\hat{\mP}$. Moreover, $\hat{\mu}^1_t=\hat{\mu}^2_t,$ a.s.$\hat{\mP}$
since $\hat{\cC}(\mR^n)$ separates bounded positive measures on $\mathscr{B}(\mR^n)$. Thus, the c\`adl\`ag property of
$\hat{\mu}^1_t, \hat{\mu}^2_t$ in $t$ admits us to get pathwise uniqueness, i.e.
$$
\hat{\mu}^1_t=\hat{\mu}^2_t, ~ t\in[0,T], ~ a.s.\hat{\mP}.
$$

Now, we show the claim. Set
\ce
\mathscr{S}:=\Big\{(\varphi, \psi)\in span\{h\phi; h\in\cC^1_b([0,T]), \phi\in\cB(\mR^{2n})\}\times\cB(\mR^{2n}); \\ \delta_t\bar{\hat{\mu}}^{11}_t(\varphi)=\delta_0\bar{\hat{\mu}}^{11}_0(\varphi)+\int_0^t\delta_s\bar{\hat{\mu}}^{11}_s(\psi-\bar{\gamma}(\cdot,\cdot)\varphi)\dif s, t\in[0,T]\Big\},
\de
and then by the dominated convergence theorem $\mathscr{S}$ is b.p.-closed. Define
\ce
&&\mathscr{G}:=span\{h\phi; h\in\cC^1_b([0,T]), \phi\in\cG\}, \\
&&(\bar{\mathscr{L}}\varphi)(t,x):=\phi(x)h'(t)+h(t)\bar{\mathscr{L}}_t\phi(x), \quad \varphi(t,x)=h(t)\phi(x)\in\mathscr{G},
\de
and then by (\ref{tenpro}) and the proof of Theorem \ref{unmeev}, we know that $\{(\varphi,\bar{\mathscr{L}}\varphi),
\varphi\in\mathscr{G}\}\subset\mathscr{S}$ and furthermore the b.p.-closure of $\{(\varphi,\bar{\mathscr{L}}\varphi),
\varphi\in\mathscr{G}\}\subset\mathscr{S}$. Moreover, Lemma 4.3 in \cite{lh} admits us to have
\ce
&&\left\{(\varphi,\bar{\mathscr{L}}\varphi),\varphi\in span\{h\phi; h\in\cC^1_b([0,T]), \phi\in \cC_c^\infty(\mR^{2n})\}\right\}\\
&\subset& the ~ b.p.-closure~ of ~ \{(\varphi,\bar{\mathscr{L}}\varphi),\varphi\in\mathscr{G}\}.
\de
Thus,
\be
\left\{(\varphi,\bar{\mathscr{L}}\varphi), \varphi\in span\{h\phi; h\in\cC^1_b([0,T]), \phi\in \cC_c^\infty(\mR^{2n})\}\right\}\subset\mathscr{S}.
\label{opebel}
\ee
Besides, it follows from Problem 4.11.12 of \cite{ek} that
\ce
&&\{(h,h'), h\in\cC^1_b([0,T])\}\\
&\subset& the ~ b.p.-closure~ of~ \left\{(\varphi,\bar{\mathscr{L}}\varphi), \varphi\in span\{h\phi; h\in\cC^1_b([0,T]), \phi\in \cC_c^\infty(\mR^{2n})\}\right\}.
\de
So,
\be
\{(h,h'), h\in\cC^1_b([0,T])\}\subset\mathscr{S}.
\label{10bel}
\ee
Combining (\ref{opebel}) with (\ref{10bel}), we obtain that
$$
\Big\{(\varphi,\bar{\mathscr{L}}\varphi), \varphi\in span\big\{h\phi; h\in\cC^1_b([0,T]), \phi\in span\{1, \cC_c^\infty(\mR^{2n})\}\big\}\Big\}\subset\mathscr{S},
$$
and i.e.
\ce
\delta_t\bar{\hat{\mu}}^{11}_t(\varphi)=\delta_0\bar{\hat{\mu}}^{11}_0(\varphi)+\int_0^t\delta_s\bar{\hat{\mu}}^{11}_s((\bar{\mathscr{L}}-\bar{\gamma}(\cdot,\cdot))\varphi)\dif s,\\ \forall\varphi\in span\{h\phi; h\in\cC^1_b([0,T]), \phi\in span\{1, \cC_c^\infty(\mR^{2n})\}\}.
\de
Again by the proof of Theorem \ref{unmeev}, it holds that
\ce
\bar{\hat{\mu}}^{11}_t(\phi)=\bar{\hat{\mu}}^{11}_0(\phi)+\int_0^t\bar{\hat{\mu}}^{11}_s((\bar{\mathscr{L}}_s-\bar{\gamma}(s,\cdot))\phi)\dif s, \quad \phi\in span\{1, \cC_c^\infty(\mR^{2n})\}.
\de
Thus, the claim is proved.

\bl\label{gj}
Suppose that $\{(\hat{\Omega}, \hat{\mathscr{F}}, \{\hat{\mathscr{F}}_t\}_{t\in[0,T]},\hat{\mP}), (\hat{\mu}_t,
\hat{W}_t, \hat{N}(\dif t, \dif u))\}$ is a weak solution of the Zakai equation. Then for $p>1$,
$$
\hat{\mE}\left[\sup\limits_{0\leq s\leq T}|\hat{\mu}_s(1)|^p\right]<\infty.
$$
\el
\begin{proof}
Since $\{(\hat{\Omega}, \hat{\mathscr{F}}, \{\hat{\mathscr{F}}_t\}_{t\in[0,T]},\hat{\mP}), (\hat{\mu}_t,
\hat{W}_t, \hat{N}(\dif t, \dif u))\}$ is a weak solution of the Zakai equation, $\hat{\mu}_t(1)$ satisfies
Eq.(\ref{zakaieq2}) for $F=1$, i.e.
\be
\hat{\mu}_t(1)&=&\hat{\mu}_0(1)+\int_0^t\hat{\mu}_s(\cL_s 1)\dif s+\int_0^t\hat{\mu}_s\left(\left(\sigma_2^{-1}(s)b_2(s,\cdot)\right)^i\right)\dif\hat{W}^i_s\no\\
&&+\int_0^t\int_{\mU_0}\hat{\mu}_{s-}\left(\lambda(s,\cdot,u)-1\right)\tilde{\hat{N}}(\dif s, \dif u)\no\\
&=&1+\int_0^t\hat{\mu}_s\left(\left(\sigma_2^{-1}(s)b_2(s,\cdot)\right)^i\right)\dif\hat{W}^i_s\no\\
&&+\int_0^t\int_{\mU_0}\hat{\mu}_{s-}\left(\lambda(s,\cdot,u)-1\right)\tilde{\hat{N}}(\dif s, \dif u).
\label{mi}
\ee
And then $\hat{\mu}_t(1)$ is a martingale. By maximal values inequality it holds that
$$
\hat{\mE}\left[\sup\limits_{0\leq s\leq T}|\hat{\mu}_s(1)|^p\right]\leq\left(\frac{p}{p-1}\right)^p\hat{\mE}[|\hat{\mu}_T(1)|^p].
$$

Next, we rewrite Eq.(\ref{mi}) as
\ce
\hat{\mu}_t(1)&=&1+\int_0^t\hat{\mu}_s(1)\cdot\frac{\hat{\mu}_s\left(\left(\sigma_2^{-1}(s)b_2(s,\cdot)\right)^i\right)}{\hat{\mu}_s(1)}\dif\hat{W}^i_s\\
&&+\int_0^t\int_{\mU_0}\hat{\mu}_{s-}(1)\cdot\frac{\hat{\mu}_{s-}\left(\lambda(s,\cdot,u)-1\right)}{\hat{\mu}_{s-}(1)}\tilde{\hat{N}}(\dif s, \dif u).
\de
Set
\ce
\hat{M}_t:=\int_0^t\frac{\hat{\mu}_s\left(\left(\sigma_2^{-1}(s)b_2(s,\cdot)\right)^i\right)}{\hat{\mu}_s(1)}\dif\hat{W}^i_s+\int_0^t\int_{\mU_0}\frac{\hat{\mu}_{s-}\left(\lambda(s,\cdot,u)-1\right)}{\hat{\mu}_{s-}(1)}\tilde{\hat{N}}(\dif s, \dif u),
\de
and then
\ce
\hat{\mu}_t(1)&=&\exp\left\{\hat{M}_t-\frac{1}{2}[\hat{M}^c, \hat{M}^c]_t\right\}\times\prod_{0<s\leq t}(1+\triangle\hat{M}_s)e^{-\triangle\hat{M}_s},
\de
where $\hat{M}^c$ denotes the continuous part of $\hat{M}$ and $\triangle\hat{M}_s=\hat{M}_s-\hat{M}_{s-}$. So, the H\"older
inequality admits us to obtain that for $p<\alpha<\frac{1}{1-l}$,
\ce
\hat{\mE}[|\hat{\mu}_T(1)|^p]&=&\hat{\mE}\left(\exp\left\{p\hat{M}_T-\frac{p}{2}[\hat{M}^c, \hat{M}^c]_T\right\}\times\prod_{0<s\leq T}(1+\triangle\hat{M}_s)^pe^{-p\triangle\hat{M}_s}\right)\\
&=&\hat{\mE}\Bigg(\left[\exp\left\{\alpha\hat{M}_T-\frac{1}{2}[\alpha\hat{M}^c, \alpha\hat{M}^c]_T\right\}\times\prod_{0<s\leq T}(1+\alpha\triangle\hat{M}_s)e^{-\alpha\triangle\hat{M}_s}\right]^{\frac{p}{\alpha}}\\
&&\times\exp\left\{\frac{(\alpha-1)p}{2}[\hat{M}^c, \hat{M}^c]_T\right\}\times\prod_{0<s\leq T}\frac{(1+\triangle\hat{M}_s)^p}{(1+\alpha\triangle\hat{M}_s)^{\frac{p}{\alpha}}}\Bigg)\\
&\leq&\left(\hat{\mE}\cE(\alpha\hat{M})_T\right)^{\frac{p}{\alpha}}\\
&&\cdot\left(\hat{\mE}\exp\left\{\frac{(\alpha-1)p}{2}\cdot\frac{\alpha}{\alpha-p}[\hat{M}^c, \hat{M}^c]_T\right\}\times\prod_{0<s\leq T}\left(\frac{(1+\triangle\hat{M}_s)^p}{(1+\alpha\triangle\hat{M}_s)^{\frac{p}{\alpha}}}\right)^{\frac{\alpha}{\alpha-p}}\right)^{\frac{\alpha-p}{\alpha}}.
\de

Let us calculate the right hand side of the above inequality. Since $\alpha\triangle\hat{M}_s>-1$, and by $(\mathbf{H}^1_{b_2, \sigma_2, f_2})$,
\ce
&&\hat{\mE}\Big[\exp\Big\{\frac{1}{2}[\hat{M}^c,\hat{M}^c]_T+<\hat{M}^d,\hat{M}^d>_T\Big\}\Big]\\
&=&\hat{\mE}\Big[\exp\Big\{\frac{1}{2}\int_0^T\left|\frac{\hat{\mu}_s\left(\left(\sigma_2^{-1}(s)b_2(s,\cdot)\right)^i\right)}{\hat{\mu}_s(1)}\right|^2\dif s\\
&&\qquad +\int_0^T\int_{\mU_0}\left(\frac{\hat{\mu}_{s-}\left(\lambda(s,\cdot,u)-1\right)}{\hat{\mu}_{s-}(1)}\right)^2\nu(\dif u)\dif s\Big\}\Big]\\
&\leq&\exp\Big\{\frac{1}{2}L_2^4T+T\int_{\mU_0}(1-L(u))^2\nu(\dif u)\Big\}<\infty,
\de
where $\hat{M}^d$ denotes the discontinuous part of $\hat{M}$, it follows from \cite[Theorem 6]{ppks} that $\cE(\alpha\hat{M})_t$ is an exponential martingale and $\hat{\mE}\cE(\alpha\hat{M})_T=1$. Besides, we have by $(\mathbf{H}^1_{b_2, \sigma_2, f_2})$
\ce
&&\hat{\mE}\exp\left\{\frac{(\alpha-1)p}{2}\cdot\frac{\alpha}{\alpha-p}[\hat{M}^c, \hat{M}^c]_T\right\}\times\prod_{0<s\leq T}\left(\frac{(1+\triangle\hat{M}_s)^p}{(1+\alpha\triangle\hat{M}_s)^{\frac{p}{\alpha}}}\right)^{\frac{\alpha}{\alpha-p}}\\
&=&\hat{\mE}\exp\left\{\frac{(\alpha-1)p\alpha}{2(\alpha-p)}\int_0^T\left|\frac{\hat{\mu}_s\left(\left(\sigma_2^{-1}(s)b_2(s,\cdot)\right)^i\right)}{\hat{\mu}_s(1)}\right|^2\dif s\right\}\\
&&\qquad \times\prod_{0<s\leq T}\left(\frac{(1+\triangle\hat{M}_s)^p}{(1+\alpha\triangle\hat{M}_s)^{\frac{p}{\alpha}}}\right)^{\frac{\alpha}{\alpha-p}}\\
&\leq&\exp\left\{\frac{(\alpha-1)p\alpha}{2(\alpha-p)}L_2^4T\right\}\times\hat{\mE}\prod_{0<s\leq T}\left(\frac{(1+\triangle\hat{M}_s)^p}{(1+\alpha\triangle\hat{M}_s)^{\frac{p}{\alpha}}}\right)^{\frac{\alpha}{\alpha-p}}.
\de
Set
$$
G(x):=\frac{(1+x)^p}{(1+\alpha x)^{\frac{p}{\alpha}}}, \quad l-1<x<0,
$$
and then
$$
\lim_{x\uparrow0}\frac{\log G(x)}{x^2}=\frac{p(\alpha-1)}{2}.
$$
Thus, there exists a constant $C_{\alpha, p}$ such that
$$
|\log G(x)|\leq C_{\alpha, p}x^2, \quad l-1<x<0.
$$
Based on this, one can obtain that
\ce
\hat{\mE}\left(\prod_{0<s\leq T}G(\triangle\hat{M}_s)^{\frac{\alpha}{\alpha-p}}\right)&=&\hat{\mE}\left(\exp\left\{\sum_{0<s\leq T}\frac{\alpha}{\alpha-p}\log G(\triangle\hat{M}_s)\right\}\right)\\
&\leq&\hat{\mE}\left(\exp\left\{\sum_{0<s\leq T}\frac{\alpha C_{\alpha, p}}{\alpha-p}|\triangle\hat{M}_s|^2\right\}\right)\\
&\leq&\hat{\mE}\left(\exp\left\{\frac{\alpha C_{\alpha, p}}{\alpha-p}\int_0^T\int_{\mU_0}(1-L(u))^2\hat{N}(\dif s, \dif u)\right\}\right)\\
&=&\exp\left\{\frac{\alpha C_{\alpha, p}}{\alpha-p}T\int_{\mU_0}\left(e^{(1-L(u))^2}-1\right)\nu(\dif u)\right\}<\infty,
\de
where Lemma A.2 in \cite{hqxz} is used in the last second equality. The proof is completed.
\end{proof}

{\bf Proof of Theorem \ref{unipalazakai} (ii)}: Take two weak solutions $\{(\hat{\Omega}^1, \hat{\mathscr{F}}^1, \{\hat{\mathscr{F}}^1_t\}_{t\in[0,T]}, \hat{\mP}^1), (\hat{\mu}^1_t,\hat{W}^1_t, \\\hat{N}^1(\dif t, \dif u))\}$ and $\{(\hat{\Omega}^2, \hat{\mathscr{F}}^2, \{\hat{\mathscr{F}}^2_t\}_{t\in[0,T]}, \hat{\mP}^2), (\hat{\mu}^2_t, \hat{W}^2_t, \hat{N}^2(\dif t, \dif u))\}$ with $\hat{\mP}^1\circ(\hat{\mu}^1_0)^{-1}=\hat{\mP}^2\circ(\hat{\mu}^2_0)^{-1}$. Define
\ce
&&\zeta_t^1:=\int_0^t\int_{\mU\setminus\mU_0}u\hat{N}^1(\dif s, \dif u)+\int_0^t\int_{\mU_0}u\tilde{\hat{N}}^1(\dif s, \dif u),\\
&&\zeta_t^2:=\int_0^t\int_{\mU\setminus\mU_0}u\hat{N}^2(\dif s, \dif u)+\int_0^t\int_{\mU_0}u\tilde{\hat{N}}^2(\dif s, \dif u),
\de
and then $\zeta^1, \zeta^2$ are two $\mU-$valued processes. Moreover, $\hat{\mP}^1\circ (\hat{W}^1,\zeta^1)^{-1}=\hat{\mP}^2\circ (\hat{W}^2,\zeta^2)^{-1}=:\mQ$.
Let $\hat{\mP}^1_{\hat{W}^1, \zeta^1}, \hat{\mP}^2_{\hat{W}^2, \zeta^2}$ denote conditional distribution of $\hat{\mu}^1, \hat{\mu}^2$ with respect to $(\hat{W}^1, \zeta^1), (\hat{W}^2, \zeta^2)$, respectively. Set
\ce
&&\hat{\hat{\Omega}}:=D([0,T],\cM(\mR^n))\times D([0,T],\cM(\mR^n))\times C([0,T],\mR^m)\times D([0,T],\mU),\\
&&\hat{\hat{\mP}}(\dif\omega^1,\dif\omega^2,\dif\omega^3,\dif\omega^4):=\hat{\mP}^1_{\hat{W}^1, \zeta^1}(\dif\omega^1)\hat{\mP}^2_{\hat{W}^2, \zeta^2}(\dif\omega^2)Q(\dif\omega^3,\dif\omega^4), \quad \omega=(\omega^1,\omega^2,\omega^3,\omega^4)\in\hat{\hat{\Omega}},\\
&&\hat{\hat{\cF}}:=\overline{\mathscr{B}(D([0,T],\cM(\mR^n)))\times \mathscr{B}(D([0,T],\cM(\mR^n)))\times \mathscr{B}(C([0,T],\mR^m))\times \mathscr{B}(D([0,T],\mU))}^{\hat{\hat{\mP}}},
\de
and then $(\hat{\hat{\Omega}}, \hat{\hat{\cF}}, \hat{\hat{\mP}})$ is a complete probability space. Furthermore, put
\ce
&&\mathscr{B}_t(D([0,T],\cM(\mR^n))):=\sigma(\omega^1(s), 0\leq s\leq t, \omega^1\in D([0,T],\cM(\mR^n))),\\
&&\mathscr{B}_t:=\mathscr{B}_t(D([0,T],\cM(\mR^n)))\times \mathscr{B}_t(D([0,T],\cM(\mR^n)))\times \mathscr{B}_t(C([0,T],\mR^m))\times \mathscr{B}_t(D([0,T],\mU)),\\
&&\hat{\hat{\cF}}_t:=\bigcap\limits_{\varepsilon>0}\sigma(\mathscr{B}_{t+\varepsilon}, \cN),
\de
where $\cN:=\{A\subset\hat{\hat{\Omega}}| A\subset B, \hat{\hat{\mP}}(B)=0\}$, and then $(\hat{\hat{\Omega}}, \hat{\hat{\cF}}, \hat{\hat{\mP}}, \{\hat{\hat{\cF}}_t\}_{t\in[0,T]})$ is a complete filtered probability space. On the space we have

(i)
$$
\hat{\hat{\mP}}\circ (\omega^1,\omega^3,\omega^4)^{-1}=\hat{\mP}^{1}\circ (\hat{\mu}^1,\hat{W}^1,\zeta^1)^{-1},
$$

(ii)
$$
\hat{\hat{\mP}}\circ (\omega^2,\omega^3,\omega^4)^{-1}=\hat{\mP}^{2}\circ (\hat{\mu}^2,\hat{W}^2,\zeta^2)^{-1},
$$

(iii) $\omega^3$ is a Brownian motion,

(iv) $\hat{\hat{N}}(\dif t, \dif u)$ is a Poisson random measure with the predictable compensator $\dif t\nu(\dif u)$,
where $\hat{\hat{N}}((0,t], \dif u):=\sharp\{0<s\leq t: \omega^4(s)-\omega^4(s-)\in\dif u\}$ and ``$\sharp$" stands for cardinal number.

Next, combining (i), (ii), (iii) with (iv), we obtain that $\{(\hat{\hat{\Omega}}, \hat{\hat{\mathscr{F}}}, \{\hat{\hat{\mathscr{F}}}_t\}_{t\in[0,T]}, \hat{\hat{\mP}}), (\omega^1,\omega^3, \\\hat{\hat{N}}(\dif t, \dif u))\}$ and $\{(\hat{\hat{\Omega}}, \hat{\hat{\mathscr{F}}}, \{\hat{\hat{\mathscr{F}}}_t\}_{t\in[0,T]}, \hat{\hat{\mP}}), (\omega^2,\omega^3, \hat{\hat{N}}(\dif t, \dif u))\}$ are two weak solutions. It follows from Theorem \ref{unipalazakai} (i) that
\ce
\hat{\hat{\mP}}\{\omega^1=\omega^2\}=1.
\de
Thus, for any $\Gamma\in\mathscr{B}(D([0,T],\cM(\mR^n)))\times \mathscr{B}(C([0,T],\mR^m))\times \mathscr{B}(D([0,T],\mU))$,
\ce
\hat{\mP}^{1}\{(\hat{\mu}^1,\hat{W}^1,\zeta^1)\in\Gamma\}&=&\hat{\hat{\mP}}\{(\omega^1,\omega^3,\omega^4)\in\Gamma\}=\hat{\hat{\mP}}\{(\omega^2,\omega^3,\omega^4)\in\Gamma\}\\
&=&\hat{\mP}^{2}\{(\hat{\mu}^2,\hat{W}^2,\zeta^2)\in\Gamma\}.
\de
That is, $\{(\hat{\mu}^1_t,\hat{W}^1_t, \hat{N}^1(\dif t, \dif u)), t\in[0,T]\}$ and $\{(\hat{\mu}^2_t,\hat{W}^2_t, \hat{N}^2(\dif t, \dif u)), t\in[0,T]\}$
have the same finite-dimensional distributions.

\section{Uniqueness in joint law for the Kushner-Stratonovich equation} \label{uniqueks}

In the section, these definitions of weak solution and uniqueness in joint law for the Kushner-Stratonovich equation are given. Later, we
prove that weak solutions for the Kushner-Stratonovich equation are unique in joint law.

\bd\label{soluks}
If there exists the pair $\{(\bar{\Omega}, \bar{\mathscr{F}}, \{\bar{\mathscr{F}}_t\}_{t\in[0,T]}, \bar{\mP}), (\pi_t,
I_t, U(\dif t, \dif u))\}$ such that the following hold:

(i) $(\bar{\Omega}, \bar{\mathscr{F}}, \{\bar{\mathscr{F}}_t\}_{t\in[0,T]},\bar{\mP})$ is a complete filtered
probability space;

(ii) $\pi_t$ is a $\cP(\mR^n)$-valued $\bar{\mathscr{F}}_t$-adapted c\`adl\`ag process;

(iii) $I_t$ is a $m$-dimensional $\bar{\mathscr{F}}_t$-adapted Brownian motion;

(iv) $U(\dif t, \dif u)$ is a Poisson random measure with a predictable compensator \\
$\pi_{t}\left(\lambda(t,\cdot,u)\right)\dif t\nu(\dif u)$;

(v) $(\pi_t, I_t, U(\dif t, \dif u))$ satisfies the following equation
\be
\pi_t(F)&=&\pi_0(F)+\int_0^t\pi_s(\cL_s F)\dif s+\int_0^t\bigg(\pi_s\left(F\left(\sigma_2^{-1}(s)b_2(s,\cdot)\right)^i\right)\no\\
&&\qquad\qquad -\pi_s\left(F\right)\pi_s\left(\left(\sigma_2^{-1}(s)b_2(s,\cdot)\right)^i\right)\bigg)\dif I^i_s\no\\
&&+\int_0^t\int_{\mU_0}\frac{\pi_{s-}\left(F\lambda(s,\cdot,u)\right)-\pi_{s-}\left(F\right)\pi_{s-}
\left(\lambda(s,\cdot,u)\right)}{\pi_{s-}\left(\lambda(s,\cdot,u)\right)}\tilde{U}(\dif s, \dif u),\no\\
&&\qquad\qquad\qquad\qquad\qquad\qquad F\in \cC_c^\infty(\mR^n),
\label{ks}
\ee
where $\tilde{U}(\dif t, \dif u)=U(\dif t, \dif u)-\pi_{t}\left(\lambda(t,\cdot,u)\right)\dif t\nu(\dif u)$,
then $\{(\bar{\Omega}, \bar{\mathscr{F}}, \{\bar{\mathscr{F}}_t\}_{t\in[0,T]}, \bar{\mP}), (\pi_t,
I_t, \\U(\dif t, \dif u))\}$ is called a weak solution of the Kushner-Stratonovich equation.
\ed

By the deduction in Subsection \ref{nonfilter}, it is obvious that $\{(\Omega, \mathscr{F}, \{\mathscr{F}_t\}_{t\in[0,T]}, \mP),
(\mP_t, \bar{W}_t,\\ N_\lambda(\dif t, \dif u))\}$ is a weak solution of the Kushner-Stratonovich equation.

\bd\label{launks}
Uniqueness in joint law of the Kushner-Stratonovich equation means that if there exist two weak solutions $\{(\bar{\Omega}^1, \bar{\mathscr{F}}^1, \{\bar{\mathscr{F}}^1_t\}_{t\in[0,T]}, \bar{\mP}^1), (\pi^1_t, I^1_t, U^1(\dif t, \dif u))\}$ and $\{(\bar{\Omega}^2, \bar{\mathscr{F}}^2, \{\bar{\mathscr{F}}^2_t\}_{t\in[0,T]}, \bar{\mP}^2), (\pi^2_t, I^2_t, U^2(\dif t, \dif u))\}$ with $\bar{\mP}^1\circ(\pi^1_0)^{-1}=\bar{\mP}^2\circ(\pi^2_0)^{-1}$, then \\
$\{(\pi^1_t, I^1_t, U^1(\dif t, \dif u)), t\in[0,T]\}$ and $\{(\pi^2_t, I^2_t, U^2(\dif t, \dif u)), t\in[0,T]\}$
have the same finite-dimensional distributions.
\ed

Here, we give out the main result in the section.

\bt\label{ksun}
Suppose that {\bf Assumption 1.-2.} are satisfied. Then the Kushner-Stratonovich equation has uniqueness in joint law.
\et
\begin{proof}
Take two weak solutions $\{(\bar{\Omega}^1, \bar{\mathscr{F}}^1, \{\bar{\mathscr{F}}^1_t\}_{t\in[0,T]}, \bar{\mP}^1), (\pi^1_t,
I^1_t, U^1(\dif t, \dif u))\}$ and $\{(\bar{\Omega}^2, \bar{\mathscr{F}}^2,\\ \{\bar{\mathscr{F}}^2_t\}_{t\in[0,T]},\bar{\mP}^2), (\pi^2_t,
I^2_t, U^2(\dif t, \dif u))\}$ with $\bar{\mP}^1\circ(\pi^1_0)^{-1}=\bar{\mP}^2\circ(\pi^2_0)^{-1}$. And then set
\ce
\frac{1}{\chi_t^1}&:=&\exp\bigg\{-\int_0^t\pi^1_s\left(\left(\sigma_2^{-1}(s)b_2(s,\cdot)\right)^i\right)\dif I^{1,i}_s-\frac{1}{2}\sum\limits_{i=1}^m\int_0^t
\left|\pi^1_s\left(\left(\sigma_2^{-1}(s)b_2(s,\cdot)\right)^i\right)\right|^2\dif s\\
&&\quad\qquad -\int_0^t\int_{\mU_0}\log\pi^1_{s-}\left(\lambda(s,\cdot,u)\right)\tilde{U}^1(\dif s, \dif u)\\
&&\quad\qquad-\int_0^t\int_{\mU_0}\left[\log\pi^1_{s}\left(\lambda(s,\cdot,u)\right)+\frac{(1-\pi^1_{s}\left(\lambda(s,\cdot,u)\right))}{\pi^1_{s}\left(\lambda(s,\cdot,u)\right)}\right]\pi^1_{s}\left(\lambda(s,\cdot,u)\right)\nu(\dif u)\dif s\bigg\},
\de
and by the similar deduction to $\Lambda_t$ in Subsection \ref{nonfilter}, we know that $\frac{1}{\chi_t^1}$ is an exponential martingale. Define a probability measure $\mQ^1$ via
$$
\frac{\dif \mQ^1}{\dif \bar{\mP}^1}=\frac{1}{\chi_T^1}.
$$
By the Girsanov theorem for Brownian motions and random measures, it holds that under the measure $\mQ^1$
$$
\hat{W}^{1,i}_t:=I^{1,i}_t+\int_0^t\pi^1_s\left(\left(\sigma_2^{-1}(s)b_2(s,\cdot)\right)^i\right)\dif s
$$
is a Brownian motion and
$$
\tilde{\bar{U}}^1(\dif t, \dif u):=U^1(\dif t, \dif u)-\nu(\dif u)\dif t
$$
is a Poisson compensated martingale measure. Note that for any $F\in \cC_c^\infty(\mR^n)$, $\pi^1_t(F)$ satisfies Eq.(\ref{ks}).
Thus, applying the It\^o formula to $\pi^1_t(F)\chi_t^1$, we obtain that
\ce
\pi^1_t(F)\chi_t^1&=&\pi^1_0(F)+\int_0^t\pi^1_s(\cL_s F)\chi_s^1\dif s+\int_0^t\pi^1_s\left(F\left(\sigma_2^{-1}(s)b_2(s,\cdot)\right)^i\right)\chi_s^1\dif\hat{W}^{1,i}_s\no\\
&&+\int_0^t\int_{\mU_0}\pi^1_s\(F(\lambda(s, \cdot, u)-1)\)\chi_s^1\tilde{\bar{U}}^1(\dif s, \dif u).
\de
Set
$$
\hat{\mu}^1_t:=\pi^1_t\chi_t^1,
$$
and then $\{(\bar{\Omega}^1, \bar{\mathscr{F}}^1, \{\bar{\mathscr{F}}^1_t\}_{t\in[0,T]}, \mQ^1), (\hat{\mu}^1_t,
\hat{W}^1_t, U^1(\dif t, \dif u))\}$ is a weak solution of the Zakai equation. By the same way to above, we can
define $\chi_t^2, \mQ^2, \hat{\mu}^2_t, \hat{W}^2_t, \tilde{\bar{U}}^2(\dif t, \dif u)$ and furthermore attain that
$\{(\bar{\Omega}^2, \bar{\mathscr{F}}^2, \{\bar{\mathscr{F}}^2_t\}_{t\in[0,T]}, \mQ^2), (\hat{\mu}^2_t,
\hat{W}^2_t, U^2(\dif t, \dif u))\}$ is also a weak solution of the Zakai equation. Thus, it follows from Theorem \ref{unipalazakai} (ii)
that $(\hat{\mu}^1, \hat{W}^1, U^1)$ and $(\hat{\mu}^2, \hat{W}^2, U^2)$ have the same finite dimensional distributions. 
Moreover, define
\ce
&&\bar{\eta}_t^1:=\int_0^t\int_{\mU\setminus\mU_0}uU^1(\dif s, \dif u)+\int_0^t\int_{\mU_0}u\tilde{\bar{U}}^1(\dif s, \dif u),\\
&&\bar{\eta}_t^2:=\int_0^t\int_{\mU\setminus\mU_0}uU^2(\dif s, \dif u)+\int_0^t\int_{\mU_0}u\tilde{\bar{U}}^2(\dif s, \dif u),
\de
and then $(\hat{\mu}^1, \hat{W}^1, \bar{\eta}^1)$ and $(\hat{\mu}^2, \hat{W}^2, \bar{\eta}^2)$ also have the same finite dimensional distributions.  

Next, by the deduction in Subsection \ref{nonfilter}, one know that
\ce
\pi^1_t(F)=\frac{\hat{\mu}^1_t(F)}{\hat{\mu}^1_t(1)}, \qquad F\in \cC_c^\infty(\mR^n),
\de
and furthermore
$$
I^{1,i}_t=\hat{W}^{1,i}_t-\int_0^t\frac{\hat{\mu}^1_s\left(\left(\sigma_2^{-1}(s)b_2(s,\cdot)\right)^i\right)}{\hat{\mu}^1_s(1)}\dif s.
$$
Define
\ce
&&\eta_t^1:=\int_0^t\int_{\mU\setminus\mU_0}uU^1(\dif s, \dif u)+\int_0^t\int_{\mU_0}u\tilde{U}^1(\dif s, \dif u),\\
&&\eta_t^2:=\int_0^t\int_{\mU\setminus\mU_0}uU^2(\dif s, \dif u)+\int_0^t\int_{\mU_0}u\tilde{U}^2(\dif s, \dif u),
\de
and then $\eta^1, \eta^2$ are two $\mU-$valued processes under $\bar{\mP}^1, \bar{\mP}^2$, respectively. Moreover, under $\mQ^1, \mQ^2$, $\eta^1, \eta^2$ become $\bar{\eta}^1, \bar{\eta}^2$, respectively. Thus, there exists a measurable functional $\Psi: D([0,T],\cM(\mR^n))\times C([0,T],\mR^m)\times D([0,T],\mU)
\mapsto D([0,T],\cP(\mR^n))\times C([0,T],\mR^m)\times D([0,T],\mU)$ such that
$$
(\pi^1, I^{1}, \eta^1)=\Psi(\hat{\mu}^1, \hat{W}^{1}, \bar{\eta}^1).
$$
By the same deduction to above, it holds that
$$
(\pi^2, I^{2}, \eta^2)=\Psi(\hat{\mu}^2, \hat{W}^{2}, \bar{\eta}^2).
$$
And then for $\Gamma\in \mathscr{B}(D([0,T],\cP(\mR^n)))\times\mathscr{B}(C([0,T],\mR^m))\times\mathscr{B}(D([0,T],\mU))$,
\ce
\bar{\mP}^1((\pi^1, I^{1}, \eta^1)\in\Gamma)&=&\mE_{\mQ^1}[I_{\{\Psi(\hat{\mu}^1, \hat{W}^{1}, \bar{\eta}^1)\in \Gamma\}}\hat{\mu}^1_T(1)]
=\mE_{\mQ^2}[I_{\{\Psi(\hat{\mu}^2, \hat{W}^{2}, \bar{\eta}^2)\in \Gamma\}}\hat{\mu}^2_T(1)]\\
&=&\bar{\mP}^2((\pi^2, I^{2}, \eta^2)\in\Gamma),
\de
where $\mE_{\mQ^1}$ stands for the expectation under the probability $\mQ^1$. So, $\{(\pi^1_t, I^1_t, U^1(\dif t, \dif u)), t\in[0,T]\}$ and $\{(\pi^2_t, I^2_t, U^2(\dif t, \dif u)), t\in[0,T]\}$ have the same finite-dimensional distributions.
\end{proof}

\br
Here we can't obtain pathwise uniqueness of the Kushner-Stratonovich equation from uniqueness in joint law for weak solutions, since the innovations conjecture analogue doesn't hold. (See \cite[Proposition 5.9]{lh2} for details)
\er

\bigskip

\textbf{Acknowledgements:}

The author would like to thank Professor Xicheng Zhang for his valuable discussions. And the author would also
wish to thank two anonymous referees for giving useful suggestions to improve this paper.

\end{document}